\documentclass[twoside,11pt,reqno]{article}
\usepackage{amssymb}
\usepackage{amsthm}
\usepackage[tbtags]{amsmath}
\usepackage{doc}
\usepackage{latexsym}
\usepackage{amscd}
  \theoremstyle{change}

\theoremstyle{change}  
\begin{document}
\thispagestyle{plain}
 \markboth{}{}
\small{\addtocounter{page}{0} \pagestyle{plain}
 \vspace{0.2in}
\noindent{\large \bf On the Twisted $q$-Euler numbers and
polynomials associated with basic $q$-$l$-functions }
\vspace{0.25in}\\
\noindent{\sc Taekyun Kim,}
{\it EECS, Kyungpook National University, Daegu 702-701, S. Korea\\
e-mail} : {\verb|tkim64@hanmail.net|} ~{\it or}~
{\verb|tkim@knu.ac.kr|}
\vspace{0.15in}\\
\noindent{\sc Seog-Hoon Rim,}
{\it Department of Mathematics Education, Kyungpook National University, Daegu 702-701, S. Korea\\
e-mail} : {\verb|shrim@knu.ac.kr|}
\vspace{0.15in}\\
{\footnotesize {\sc Abstract.} One purpose of this paper is to
construct twisted $q$-Euler numbers by using $p$-adic invariant
integral on $\mathbb{Z}_{p}$ in the sense of fermionic. Finally, we
consider twisted Euler $q$-zeta function and $q$-$l$-series which
interpolate twisted $q$-Euler numbers and polynomials.}
\vspace{0.2in}\\
\noindent{\bf 1. Introduction} \setcounter{equation}{0}
\vspace{0.1in}\\
\indent Let $p$ be a fixed odd positive integer. Throughout this
paper $\mathbb{Z}_{p}$, $\mathbb{Q}_{p}$, $\mathbb{C}$ and
$\mathbb{C}_{p}$ are respectively denoted as the ring of $p$-adic
rational integers, the field of $p$-adic rational numbers, the
complex numbers field and the completion of algebraic closure of
$\mathbb{Q}_{p}$. The $p$-adic absolute value in $\mathbb{C}_{p}$
is normalized so that $|p|_{p}=1/p$. When one talks of
$q$-extension, $q$ is considered in many ways such as an
indeterminate, a complex number $q \in \mathbb{C}$, or a $p$-adic
number $q \in \mathbb{C}_{p}$. If $q \in \mathbb{C}$, one normally
assumes that $|q|<1$. If $q \in \mathbb{C}_{p}$, we normally
assume that $|q-1|_{p} < p^{-1/p-1}$ so that $q^{x} = \exp(x \log
q)$ for $|x|_{p} \leq 1$, cf. [1, 2, 15, 20].

\noindent We use the notations as
\begin{eqnarray}
\lceil x\rceil _{q} &=& \dfrac{1-q^{x}}{1-q} = 1 + q + q^{2} +
\cdots + q^{x-1}, \notag \\
\lceil x \rceil_{-q} &=& \dfrac{1-(-q)^{x}}{1+q} \notag \\
&=& 1 - q + q^{2} - q^{3} + \cdots + (-1)^{x} q^{x-1}, \quad
\text{see [12]} \notag.
\end{eqnarray}
\indent Let $UD(\mathbb{Z}_{p})$ be the set of uniformly
differentiable function on $\mathbb{Z}_{p}$. For $f \in
UD(\mathbb{Z}_{p})$, the $p$-adic invariant $q$-integral on
$\mathbb{Z}_{p}$ was defined by
$$I_{q}(f)= \int_{\mathbb{Z}_{p}} f(x) d\mu_{q}(x) = \lim_{\mathbf{v} \to
\infty} \dfrac{1}{[p^{\mathbf{v}}]_{q}}
\sum^{p^{\mathbf{v}}-1}_{x=0} f(x) q^{x}, \quad \text{(see [3, 8,
9])}.$$ Note that $$I_{1}(f) = \lim_{q \to 1} I_{q}(f) = \int
_{\mathbb{Z}_{p}} f(x) d\mu_{1}(x) = \lim_{\mathbf{v} \to \infty}
\dfrac{1}{p^{\mathbf{v}}} \sum_{x=0}^{p^{\mathbf{v}}-1} f(x).$$
\noindent For a fixed positive integer $d$ with $(p,d)=1$, set
\begin{eqnarray}
X=X_{d} &=& \lim_{\overleftarrow{N}}
\mathbb{Z}/dp^{\mathbf{v}}\mathbb{Z}
\notag \\
X_{1} &=& \mathbb{Z}_{p} \notag \\
X^{\ast} &=& \bigcup_{0 < a < dp, (a,p)=1 } \notag \\
a+dp^{\mathbf{v}} \mathbb{Z}_{p} &=& \{ x \in X ~ |~  x \equiv a
~({\rm mod} ~ dp^{\mathbf{v}})\}, \notag
\end{eqnarray}
where $a \in \mathbb{Z}$ satisfies the condition $0 \leq a <
dp^{\mathbf{v}}$, cf. [15].

Recently Kim consider $q$-Euler numbers as follows:
$$F_{q}(x,t)=\lceil 2 \rceil_{q} \sum_{n=0}^{\infty} (-1)^{n} q^{n}
e^{\lceil n+x \rceil_{q} t} = \sum_{n=0}^{\infty} E_{n,q}(x)
\dfrac{t^{n}}{n!}, \qquad \text{see [5, 12]}.$$

\noindent For $h \in \mathbb{Z}$, $k \in \mathbb{N}$, $q$-Euler
polynomials of higher order were defined by
\begin{equation}
E_{n,q}^{(h,k)}(x) = \underbrace{ \int_{\mathbb{Z}_{p}} \cdots
\int_{\mathbb{Z}_{p}}}_{k-\text{times}} \lceil x+x_{1} + \cdots +
x_{k} \rceil_{q}^{n} q^{\sum_{i=1}^{k} x_{i}(h-i)}
d\mu_{-q}(x_{1}) \cdots d\mu_{-q}(x_{k}), ~\text{see [8, 15]}.
\end{equation}
\noindent The $q$-Euler polynomial at $x=0$ are called $q$-Euler
numbers, cf. [1, 2, 14]. Hence $q$-Euler numbers of higher order
can be defined by $E_{n,q}^{(h,k)} = E_{n,q}^{(h,k)}(0)$.
\\
\indent Let $\chi$ be a Dirichlet Character with conductor
$f(=\text{odd}) \in \mathbb{N}$. Then the $l$-function attached to
$\chi$ is defined as
$$l(s,\chi) = 2 \sum_{n=1}^{\infty} \dfrac{\chi(n)(-1)^{n}}{n^{s}},
\quad ~ \text{for} ~ s \in \mathbb{C}, \quad ~ \text{see [4,
13]}.$$ \noindent When $\chi=1$, this is the Euler zeta function,
which is defined by $$\zeta_{E}(s) = 2 \sum_{n=1}^{\infty}
\dfrac{(-1)^{n}}{n^{s}}, \quad ~ \text{see [3, 4, 13]}.$$
\noindent In [8, 12], $q$-analogue of Euler $\zeta$-function is
defined by
\begin{equation}
\zeta_{E,q}^{(h)} (s,x) = \lceil 2 \rceil _{q} \sum_{n=0}^{\infty}
\dfrac{(-1)^{n}q^{nh}}{\lceil n+x \rceil _{q}^{s}}, \quad ~ s, h
\in \mathbb{C}.
\end{equation}
\noindent Note that $\zeta_{E,q}^{(h)}(s,x)$ is analytic
continuation an $\mathbb{C}$ with
\begin{equation}
\zeta_{E,q}^{(h)}(-n,x) = E_{n,q}^{(h,1)}(x), \quad ~ \text{for} ~ n
\in \mathbb{N},~ h\in \mathbb{Z}, \quad ~ \text{see [8, 5]}.
\end{equation}
\noindent It was known that
\begin{equation}
E_{n,q}^{(h,1)}(x) = \lceil 2 \rceil_{q} \sum_{k=0}^{\infty}
(-1)^{k} q^{kh} \lceil k+x \rceil_{q}^{n}, \quad ~ \text{see [8]}.
\end{equation}
\noindent It follows from (2) that $$\lim_{q \to 1}
\zeta_{E,q}^{(h)}(s,x) = \zeta_{E}(s,x) = 2 \sum_{n=0}^{\infty}
\dfrac{(-1)^n}{(n+x)^{s}}, \quad ~ \text{see [13]}.$$ \noindent
With meaning of $q$-analogue of Dirichlet's type $l$-function, we
consider the following $q$-$l$-function : for $s \in \mathbb{C}$
\begin{equation}
l_{q}(s,\chi) = \lceil 2 \rceil_{q} \sum_{k=1}^{\infty}
\dfrac{\chi(k) (-1)^{k} q^{k}}{\lceil k \rceil_{q}^{s}}, \quad ~
\text{see [3, 5, 12]}.
\end{equation}
\noindent It is easy to see that $l_{q}(s,\chi)$ is analytic
continuation in whole complex plane.

The generalized $q$-Euler numbers attached to $\chi$ were also
defined by
\begin{equation}
\lceil 2 \rceil_{q} \sum_{n=0}^{\infty} e^{\lceil n \rceil_{q}t}
\chi(n) (-1)^{n} q^{n} = \sum_{n=0}^{\infty} E_{n, \chi, q}
\dfrac{t^{n}}{n!}, \quad ~ \text{see [5, 6, 12]}.
\end{equation}
\indent Note that
\begin{equation}
\int_{\mathbb{Z}_{p}} \chi(x) \lceil x \rceil_{q}^{n} d\mu_{-q}(x) =
E_{n, \chi, q} \quad ~ \text{and}~ \quad  l_{q}(-n, x)=E_{n, \chi,
q}, \quad ~ \text{see [3, 5]}.
\end{equation}
\noindent For $m \geq 0$, extended generalized $q$-Euler numbers
attached to $\chi$ are defined as
\begin{eqnarray}
E_{n, \chi, q}^{(h, k)} &=& \underbrace{ \int_{X} \cdots \int_{X}
}_{k-\text{times}} \lceil \sum_{l=1}^k x_{l} \rceil_{q}^{n}
q^{\sum_{j=1}^{k}(h-j)x_{j}}
 (\prod_{j=1}^{k} \chi(x_{j}))d\mu_{-q}(x_{1}) \cdots d\mu_{-q}(x_{k}) \notag \\
&=& \dfrac{\lceil 2 \rceil_{q}^{k}}{\lceil 2 \rceil_{q^{f}}^{k}}
\lceil f \rceil_{q}^{n} \sum_{a_{1}, \cdots a_{k}=0}^{f-1}
q^{\sum_{l=1}^{k}(h-l+1)a_{l}} (-1)^{a_{1} + \cdots + a_{k}}
\notag \\
 & & E_{n, q^{f}}^{(h,k)} \big( \dfrac{a_l +a_2+\cdots+a_k }{f} \big) \prod_{i=1}^{k} \chi(a_{i}), ~ \quad ~ \text{(see
[8, 10]).}
\end{eqnarray}
\noindent From (8), we derive
\begin{eqnarray}
E_{n, \chi, q}^{(h,1)} &=& \dfrac{\lceil 2 \rceil_{q}}{\lceil
 2\rceil_{q^{f}}} \lceil f \rceil_{q}^{n} \sum_{a=0}^{f-1}
 q^{ha}(-1)^{a} E_{n, q^{f}}^{(h,1)}\big( \dfrac{a}{f} \big) \chi(a)
 \notag \\
&=& \lceil 2 \rceil_{q} \sum_{k=1}^{\infty} \chi(k) (-1)^{k}
q^{hk} \lceil k \rceil_{q}^{n}, \qquad ~ \text{see [5: p.5].}
\end{eqnarray}
\indent From (9), we can consider the following $q$-$l$-function:
$$l_{q}^{(h)}(s, \chi) = \lceil 2 \rceil_{q} \sum_{k=1}^{\infty}
\dfrac{\chi(k)(-1)^{k} q^{hk}}{\lceil k \rceil_{q}^{s}}, \quad ~
\text{see [5, 8, 10]},$$ where $h, s \in \mathbb{C}$. Note that
$l_{q}^{(h)}(-n,\chi) = E_{n, \chi, q}$, for $n\in \mathbb{N}, ~
h\in \mathbb{Z}$.

In the present paper, we give twisted $q$-Euler numbers by using
$p$-adic invariant $q$-integral on $\mathbb{Z}_{p}$ in the sense
of fermionic. Moreover, we construct the $q$-analogues of Euler
zeta function and $q$-$l$-function, which interpolate twisted
$q$-Euler number and polynomials at negative integer.
\vspace{0.2in}\\
\noindent{\bf 2. Twisted $q$-extension of Euler numbers}
\vspace{0.1in}\\
\indent In this section, we assume that $q \in \mathbb{C}_{p}$
with $|1-q|_{p} < 1$. \vspace{0.1in} From the definition of
$p$-adic invariant $q$-integral on $\mathbb{Z}_{p}$ in the sense
of fermionic, we derive
\begin{equation}
\lim_{q \to -1} I_{q}(f_{1}) = I_{-1}(f_{1}) = -I_{-1}(f) + 2f(0), ~
\text{where}~ f_{1}(x) = f(x+1), ~ \text{see [7]}.
\end{equation}

\indent Let $T_{p} = \cup_{n \geq 1} C_{p^{n}} = \lim_{n \to
\infty} C_{p^{n}} = C_{p^{\infty}}$ be locally constant space,
where $C_{p^{n}} = \{ w ~ | ~ w^{p^{n}}=1 \}$ is the cyclic group
of order $p^{n}$. For $w \in T_{p}$, we denote by $\phi_{w} :
\mathbb{Z}_{p} \to \mathbb{C}_{p}$ the locally constant function
$x \mapsto w^{x}$.

\noindent If we take $f(x) = \phi_{w}(x) e^{tx}$, then we have
\begin{equation}
\int_{\mathbb{Z}_{p}} e^{tx} \phi_{w}(x) d\mu_{-1}(x) =
\dfrac{2}{we^{t}+1}, \quad ~\text{see [4, 11, 13,16, 17, 18, 19]}.
\end{equation}
\noindent From (10), we derive
\begin{equation}
I_{-1}(f_{n}) = (-1)^{n}I_{-1}(f) + 2
\sum_{l=0}^{n-1}(-1)^{n-1-l}f(l), \quad ~ \text{see [6, 7]},
\end{equation}
where $f_{n}(x) = f(x+n)$. \vspace{0.1in} By (12), we easily see
that
\begin{equation}
\int_{X}e^{tx}\phi_{w}(x) \chi(x) d\mu_{-1}(x) = \dfrac{2
\sum_{i=0}^{f-1}\chi(i) \phi_{w}(i) e^{it}}{w^{f}e^{ft}+1}.
\end{equation}

\indent Now we define the analogue of Euler numbers as follows :
\begin{equation}
\dfrac{2 \sum_{i=0}^{f-1} \chi(i)\phi_{w}(i)e^{it}}{w^{f}e^{ft}+1} =
\sum_{n=0}^{\infty} E_{n, \chi, w} \dfrac{t^{n}}{n!}.
\end{equation}
\noindent From (13) and (14), we note that
\begin{equation}
\int_{X} x^{n} \phi_{w}(x) \chi(x) d\mu_{-1}(x) = E_{n, \chi, w}.
\end{equation}
\noindent In the viewpoint of (15), we consider twisted $q$-Euler
numbers using $p$-adic invariant $q$-integral on $\mathbb{Z}_{p}$
in the sense of fermionic as follows :
\begin{equation}
E_{n, w, q}^{(h,1)}(x) = \int_{\mathbb{Z}_{p}} q^{(h-1)y} w^{y}
\lceil x+y \rceil_{q}^{n} d\mu_{-q}(y).
\end{equation}
\noindent Observe that $\lim_{q \to 1} E_{n, q, w}^{(h,1)}(x) =
E_{n,w}(x)$. When $x=0$, we write $E_{n,w}^{(h)}(0,q) =
E_{n,w}^{(h)}$. Note that $E_{n, w, q}^{(h,1)}$ is the twisted
form of $E_{n, q}^{(h,1)}$, see [8].

\noindent From (16), we derive
\begin{equation}
E_{n,w,q}^{(h,1)}(x) = \dfrac{\lceil 2 \rceil_{q}}{(1-q)^{n}}
\sum_{j=0}^{n} \binom{n}{j} (-1)^{j} q^{xj} \dfrac{1}{1+q^{h+l}w}
\end{equation}

\noindent Eq.(17) is equivalent to
\begin{equation}
E_{n,w,q}^{(h,1)}(x) = \lceil 2 \rceil_{q}
\sum_{k=0}^{\infty}(-1)^{k} w^{k} q^{hk} \lceil x+k \rceil_{q}^{n},
\quad ~ h\in \mathbb{Z}, ~ n\in \mathbb{N}.
\end{equation}
\noindent By (16), we see that
$$E_{n,w,q}^{(h,1)}(x) = \dfrac{\lceil
2 \rceil_{q}}{\lceil 2 \rceil}_{q^{d}} \lceil d \rceil_{q}^{n}
\sum_{a=0}^{d-1} q^{ha} w^{a} (-1)^{a}
E_{n,q^{d}}^{(h,1)}(\dfrac{x+a}{d}),$$ where $n, ~d(=\text{odd})
\in \mathbb{N}$.

Let $\chi$ be the Dirichlet character with conductor
$f(=\text{odd}) \in \mathbf{N}$. Then we define the generalized
twisted $q$-Euler numbers as follows : For $n \geq 0$,
\begin{equation}
E_{n, w, \chi, q}^{(h,1)} = \int_{X} \chi(x) q^{(h-1)x} w^{x} \lceil
x \rceil_{q}^{n} d\mu_{-q}(x).
\end{equation}
\noindent Note that $E_{n,1,\chi,q}^{(1,1)} = E_{n, \chi,q}$, see
[3, 5, 12].

\noindent From (19), we can also derive
\begin{eqnarray}
E_{n,w,\chi,q}^{(h,1)} &=& \lceil f \rceil_{q}^{n} \dfrac{\lceil 2
\rceil_{q}}{\lceil 2 \rceil_{q^{f}}} \sum_{a=0}^{f-1} q^{ha} w^{a}
\chi(a) (-1)^{a} E_{n, w^{f}, g^{f}}^{(h,1)} \big( \dfrac{a}{f}
\big) \notag \\
&=& \lceil 2 \rceil_{q} \sum_{k=1}^{\infty} \chi(k) (-1)^{k}q^{hk}
w^{k} \lceil k \rceil_{q}^{n}.
\end{eqnarray}
\noindent It is east to check that
$$\lim_{q \to 1} E_{n, w, \chi, q}^{(h,1)} = E_{n,w,\chi}, \quad ~
\text{see [4]}.$$
\vspace{0.2in}\\
\noindent{\bf 3. $q$-Euler zeta function and $q$-$l$-functions}
\vspace{0.1in}\\
\indent In this section, we assume that $q\in \mathbb{C}$ with
$|q|<1$. Here we construct twisted $q$-Euler zeta function and
twisted $q$-$l$-function. Let $\mathbb{R}$ be the field of real
numbers and let $w$ be the $p^{r}$-th root of unity. For $q \in
\mathbb{R}$ with $0<q<1$, $s \in \mathbb{C}$ and $h\in
\mathbb{C}$, we define twisted Euler $q$-zeta function as follows
:
$$\zeta_{E,q,w}^{(h,1)}(s) = \lceil 2 \rceil_{q} \sum_{k=1}^{\infty}
\dfrac{(-1)^{k} w^{k} q^{hk}}{\lceil x+k \rceil_{q}^{s}}.$$
\noindent Note that $\zeta_{E,q,w}^{(h,1)}$ is an analytic
continuation in $\mathbb{C}$. \vspace{0.1in} In [8], it is easy to
see that
$$\zeta_{E,q,1}^{(h,1)}(s)=\zeta_{E,q}^{(h)}(s).$$
\noindent By (18), we easily see that $$\zeta_{E,q,w}^{(h,1)}(-n)
= E_{n,w,q}^{(h,1)}, \quad ~ \text{for} ~ n\in \mathbb{N}, ~ h \in
\mathbb{Z}.$$ \indent We now also consider Hurwitz's type twisted
$q$-Euler zeta function as follows : For $s \in \mathbb{C}$,
define
\begin{equation}
\zeta_{E,q,w}^{(h,1)}(s,x) = \lceil 2 \rceil_{q} \sum_{k=0}^{\infty}
\dfrac{(-1)^{k} q^{hk} w^{k}}{\lceil x+k \rceil_{q}^{s}}, \quad
~\text{for} ~ s\in \mathbb{C}, ~ h\in \mathbb{C}.
\end{equation}
\noindent Note that $\zeta_{E,q,w}^{(h,1)}(s,x)$ has an analytic
continuation on $\mathbb{C}$. By (18), we easily see that
$$\zeta_{E,q,w}^{(h,1)}(-n, x) = E_{n,w,q}^{(h,1)}(x), \quad ~
\text{for} ~ n \in \mathbb{N}, ~ h \in \mathbb{Z}.$$
\vspace{0.1in} \noindent From [8], we note that
$$\zeta_{E,q,1}^{(h,1)}(s,x) = \zeta_{E,q}^{(h)}(s,x).$$

Let $\chi$ be Dirichlet's character with conductor $f(=\text{odd})
\in \mathbb{N}$. Then we define twisted $q$-$l$-function which
interpolates twisted generalized $q$-Euler numbers attached to
$\chi$ as follows : For $s \in \mathbb{C}$, $h \in \mathbb{C}$,
define
\begin{equation}
l_{q,w}^{(h,1)}(s, \chi) = \lceil 2 \rceil_{q} \sum_{k=1}^{\infty}
\dfrac{\chi(k) (-1)^{k} q^{hk} w^{k}}{\lceil k \rceil_{q}^{s}}.
\end{equation}
\noindent For any positive integer $n$, we obtain
$$l_{q,w}^{(h,1)}(-n, \chi) = E_{n,w,\chi,q}^{(h,1)}, \quad ~
\text{for} ~ n \in \mathbb{N}, ~ h\in \mathbb{Z}.$$ \noindent From
(22), we can also derive
\begin{eqnarray}
l_{q,w}^{(h,1)}(s, \chi) &=& \lceil 2 \rceil_{q} \sum_{k=1}^{\infty}
\dfrac{\chi(k) (-1)^{k} q^{hk} w^{k}}{\lceil k \rceil_{q}^{s}}
\notag \\
&=& \lceil f \rceil_{q}^{-s} \dfrac{\lceil 2 \rceil_{q}}{\lceil 2
\rceil_{q^{f}}} \sum_{a=1}^{f} \chi(a) (-1)^{a} q^{ha} w^{a}
\zeta_{E, q^{f}, w^{f}} (s, \dfrac{a}{f}). \notag
\end{eqnarray}
{\bf Question.} Find a $q$-analogue of the $p$-adic twisted
$l$-function which interpolates generalized twisted $q$-Euler
numbers attached to $\chi$ $E_{n, w, \chi, q}^{(h,1)}$ [4, 5, 10,
13]. \vspace{0.2in} \footnotesize{

\end{document}